\documentclass{article}[10pt]

\usepackage{amssymb,latexsym,amsmath}

\usepackage{graphicx}

\begin{document}

\newcommand{\bfi}{\bfseries\itshape}

\makeatletter

\@addtoreset{figure}{section}

\def\thefigure{\thesection.\@arabic\c@figure}

\def\fps@figure{h, t}

\@addtoreset{table}{bsection}

\def\thetable{\thesection.\@arabic\c@table}

\def\fps@table{h, t}

\@addtoreset{equation}{section}

\def\theequation{\thesubsection.\arabic{equation}}

\makeatother

\newtheorem{thm}{Theorem}[section]

\newtheorem{prop}[thm]{Proposition}

\newtheorem{lema}[thm]{Lemma}

\newtheorem{cor}[thm]{Corollary}

\newtheorem{defi}[thm]{Definition}

\newtheorem{rk}[thm]{Remark}

\newtheorem{exempl}{Example}[section]

\newenvironment{exemplu}{\begin{exempl}  \em}{\hfill $\square$

\end{exempl}}

\newcommand{\comment}[1]{\par\noindent{\raggedright\texttt{#1}

\par\marginpar{\textsc{Comment}}}}

\newcommand{\todo}[1]{\vspace{5 mm}\par \noindent \marginpar{\textsc{ToDo}}\framebox{\begin{minipage}[c]{0.95 \textwidth}

\tt #1 \end{minipage}}\vspace{5 mm}\par}

\newcommand{\ea}{\mbox{{\bf a}}}

\newcommand{\eu}{\mbox{{\bf u}}}

\newcommand{\ueu}{\underline{\eu}}

\newcommand{\ueo}{\overline{u}}

\newcommand{\oeu}{\overline{\eu}}

\newcommand{\ew}{\mbox{{\bf w}}}

\newcommand{\ef}{\mbox{{\bf f}}}

\newcommand{\eF}{\mbox{{\bf F}}}

\newcommand{\eC}{\mbox{{\bf C}}}

\newcommand{\en}{\mbox{{\bf n}}}

\newcommand{\eT}{\mbox{{\bf T}}}

\newcommand{\eL}{\mbox{{\bf L}}}

\newcommand{\eR}{\mbox{{\bf R}}}

\newcommand{\eV}{\mbox{{\bf V}}}

\newcommand{\eU}{\mbox{{\bf U}}}

\newcommand{\ev}{\mbox{{\bf v}}}

\newcommand{\eve}{\mbox{{\bf e}}}

\newcommand{\uev}{\underline{\ev}}

\newcommand{\eY}{\mbox{{\bf Y}}}

\newcommand{\eK}{\mbox{{\bf K}}}

\newcommand{\eP}{\mbox{{\bf P}}}

\newcommand{\eS}{\mbox{{\bf S}}}

\newcommand{\eJ}{\mbox{{\bf J}}}

\newcommand{\eB}{\mbox{{\bf B}}}

\newcommand{\eH}{\mbox{{\bf H}}}

\newcommand{\leb}{\mathcal{ L}^{n}}

\newcommand{\eI}{\mathcal{ I}}

\newcommand{\eE}{\mathcal{ E}}

\newcommand{\hen}{\mathcal{H}^{n-1}}

\newcommand{\eBV}{\mbox{{\bf BV}}}

\newcommand{\eA}{\mbox{{\bf A}}}

\newcommand{\eSBV}{\mbox{{\bf SBV}}}

\newcommand{\eBD}{\mbox{{\bf BD}}}

\newcommand{\eSBD}{\mbox{{\bf SBD}}}

\newcommand{\ecs}{\mbox{{\bf X}}}

\newcommand{\eg}{\mbox{{\bf g}}}

\newcommand{\paromega}{\partial \Omega}

\newcommand{\gau}{\Gamma_{u}}

\newcommand{\gaf}{\Gamma_{f}}

\newcommand{\sig}{{\bf \sigma}}

\newcommand{\gac}{\Gamma_{\mbox{{\bf c}}}}

\newcommand{\deu}{\dot{\eu}}

\newcommand{\dueu}{\underline{\deu}}

\newcommand{\dev}{\dot{\ev}}

\newcommand{\duev}{\underline{\dev}}

\newcommand{\weak}{\stackrel{w}{\approx}}

\newcommand{\mild}{\stackrel{m}{\approx}}

\newcommand{\lrightarrow}{\stackrel{L}{\rightarrow}}

\newcommand{\rrightarrow}{\stackrel{R}{\rightarrow}}

\newcommand{\strong}{\stackrel{s}{\approx}}

\newcommand{\weakdown}{\rightharpoondown}

\newcommand{\opg}{\stackrel{\mathfrak{g}}{\cdot}}

\newcommand{\opunu}{\stackrel{1}{\cdot}}
\newcommand{\opdoi}{\stackrel{2}{\cdot}}

\newcommand{\opn}{\stackrel{\mathfrak{n}}{\cdot}}
\newcommand{\opx}{\stackrel{x}{\cdot}}

\newcommand{\tr}{\ \mbox{tr}}

\newcommand{\Ad}{\ \mbox{Ad}}

\newcommand{\ad}{\ \mbox{ad}}

\renewcommand{\contentsname}{ }

\title{Emergent algebras}

\author{Marius Buliga \\
\\
Institute of Mathematics of the Romanian Academy \\
P.O. BOX 1-764, RO 014700\\
Bucure\c sti, Romania\\
{\footnotesize Marius.Buliga@imar.ro}}

\date{ }

\maketitle

\begin{abstract}

Inspired from research subjects in sub-riemannian geometry and metric geometry, 
we propose uniform idempotent right quasigroups and emergent algebras as
an alternative to differentiable algebras. 

Idempotent right quasigroups (irqs) are related with racks and quandles, which appear in knot 
theory (the axioms of a irq correspond to the first two Reidemeister moves). To 
any uniform idempotent right quasigroup can be associated an approximate 
differential calculus, with Pansu differential calculus in sub-riemannian
 geometry as an example. 
 
An emergent algebra  $\mathcal{A}$ over  a uniform idempotent right 
quasigroup $X$ is a collection of operations   such that each operation   
emerges from $X$, meaning that it can be realized as a combination of the 
operations of  the uniform irq $X$, possibly by taking 
limits which are uniform with respect to a set of parameters. 

Two applications are considered: we prove a bijection between 
 contractible groups and  distributive uniform irqs (uniform quandles) and that 
 some symmetric spaces in the sense of Loos may be seen as uniform quasigroups with 
 a distributivity property.  
\end{abstract}

\noindent
{\bf Keywords:}  differential structure;  idempotent right quasigroups;
contractible groups; Carnot groups; symmetric spaces

\noindent
{\bf MSC 2000:} 20N05; 20F19; 58H15

\newpage


\section{Introduction}

A differentiable algebraic structure, or differentiable algebra, is an 
algebra (set of operations $\mathcal{A}$) over a manifold $X$ with the property that all the operations 
of the algebra are differentiable with respect to the manifold structure 
of $X$. Let us denote by $\mathcal{D}$ the differential structure of the 
manifold $X$. 
With these notations a differentiable algebra can be seen as  a pair  
$(\mathcal{A}, \mathcal{D})$, with the properties: 
\begin{enumerate}
\item[-]  $\mathcal{A}$ contains the algebraic
information, namely the collection of statements which can be formulated using 
algebraic operations, as well as  algebraic 
relations (like for example "the operation $\circ$ is associative", or 
"the operation $\circ$ is commutative", and so on),  
\item[-] $\mathcal{D}$ contains  differential structure informations, 
that is the information needed in order to 
formulate the statement "the function $f$ is differentiable", 
\item[-]  $\mathcal{A}$ 
is compatible with $\mathcal{D}$ in the sense  that any 
operation from $\mathcal{A}$ is differentiable in the sense given by 
$\mathcal{D}$. 
\end{enumerate}

In this paper we want to extend the notion of differentiable algebra, based on
the following observation: it may happen that some of the algebraic information 
encoded in $\mathcal{A}$ "emerges" from statements from 
$\mathcal{D}$. We shall give precise examples of this phenomenon in 
the paper, but our basic example is the following: the algebraic structure of a 
vector space of the tangent space at a point $x \in X$ (when $X$ is a manifold) 
is deduced from the differential structure $\mathcal{D}$. Indeed, there are 
several ways to define tangent spaces of a manifold and all are based on 
the differential structure of the manifold. For example tangent vectors may be
seen as equivalence classes of differentiable curves passing through a point, or
as derivations in the algebra of differentiable functions of the manifold.
Eventually it is proved that tangent vectors at a point of the manifold form a 
vector space. Thus from differential type information an algebraic structure 
"emerges".

 The word "emerge" is therefore  used in relation with the fact that
while  statements in $\mathcal{A}$ may be seen as constructs in universal 
algebra, statements in $\mathcal{D}$ (whatever $\mathcal{D}$ 
may exactly mean) may be obtain from the richer structure of the "differential"
world, where we can prove algebraic statements by "passing to the limit
uniformly with respect to a set of parameters". Thus, by carefully choosing 
the $\mathcal{D}$, we may obtain at least parts of $\mathcal{A}$.

We shall choose, as a minimal replacement of a  differential structure $\mathcal{D}$ on $X$, the
structure of "uniform idempotent right quasigroup", definition \ref{deftop}. 
Then a generalization of a differentiable algebra over $X$ is a algebra 
$\mathcal{A}$ over $X$ such that all operations and relations from 
$\mathcal{A}$ can be constructed or deduced from combinations of operations 
in the uniform idempotent right quasigroup $\mathcal{D}$, possibly by taking 
limits which are uniform with respect to a set of parameters. Such an algebra
may be called an "emergent algebra". In this approach, the compatibility
condition between algebraic information $\mathcal{A}$ and differential
information $\mathcal{D}$, usually expressed as the differentiability of operations 
from $\mathcal{A}$ with respect to $\mathcal{D}$, is replaced by the 
"emergence" of algebraic operations and relations from $\mathcal{D}$ (see the 
comments and references at the end of section \ref{unifirq} for a discussion 
of replacement of differentiability by emergence). 

\paragraph{Outline of the paper.} In section \ref{sec2} we explain how this emergence 
phenomenon appears in fundamental studies concerning constructions of 
the intrinsically defined  tangent bundle of a sub-riemannian manifold, like   
Bella\"{\i}che \cite{bell}, Gromov \cite{gromovsr}, Margulis and Mostow 
\cite{marmos1} \cite{marmos2}. 
In previous papers \cite{buligadil1} \cite{buligadil2} \cite{buligadil3} we
proposed a general object which allow a differential calculus on a metric 
space, called "dilatation structure", and then used it to obtain old and new
results concerning sub-riemannian geometry, in particular we used it in 
order to clarify the last section of Bella\"{\i}che paper, where he proposes to
construct the nilpotent group operation on the metric tangent space,  at a point 
of a sub-riemannian manifold, by using adapted coordinate systems and arguments 
of uniform convergence. We explain that the motivation of this paper is to 
extract the minimal algebraic (and analytical) information from dilatation 
structures, in order to  prove the appearance of emergent algebraic structures 
without using a distance (by reasoning independent of metric geometry). 

In section \ref{sec3} we introduce idempotent right quasigroups and induced operations 
(a sort of bundle of isotopes of the right quasigroup operation). In section
\ref{sec4}  
we associate to any contractible group, or to any group with a contractive
automorphism group, an idempotent right quasigroup. We show in proposition
\ref{pfirstop} that the group operation (of the contractible group) can be
reconstructed as an emergent operation from uniform limits of isotopes of 
the associated  idempotent right quasigroup. 

This result motivates us to 
introduce in section \ref{unifirq} the notion of uniform irq, definition  \ref{deftop}. The 
main property of a uniform irq is given in theorem \ref{mainthm}, namely that 
from a uniform irq $X$ we can construct a kind of tangent bundle which associates 
to any $x \in X$ a contractible group operation (like addition of vectors in the
tangent space, only that it may be non-commutative). This tangent bundle 
gives us a minimal framework for differentiability, therefore we justify the 
use of uniform irqs as a generalization of differential structures. 

In section \ref{sec6} we study two examples of emergent algebraic structures. First
example concerns contractible groups. We prove in theorem \ref{pgroudlin} 
that contractible groups are the same as distributive uniform irqs. We then 
show that such right quasigroups are in fact quasigroups. The second example
concerns symmetric spaces in the sense of Loos \cite{loos}. We prove that 
already some properties of the inverse operation in a symmetric space are true
for any uniform irq, then we propose a notion
of uniform symmetric quasigroup, seen as a uniform quasigroup with a
distributivity property for a "approximate inverse" operation $\displaystyle 
\underline{inv}_{k}$, for any $\displaystyle k \in \mathbb{N}^{\circ}$  and prove
that uniform symmetric quasigroups generate a family of symmetric spaces 
which contains the  riemannian symmetric spaces.

\section{Motivation}
\label{sec2}

A particular class of locally compact groups which admit a contractive
automorphism group is  made by Carnot groups. These groups mysteriously appear 
as models of metric tangent spaces in sub-riemannian or Carnot-Carath\'eodory 
geometry (first proved in Mitchell \cite{mit}). 

 Non-holonomic spaces were discovered in 1926 by G. Vr\u anceanu \cite{vra1},
\cite{vra2}. The  Carnot-Carath\'eodory distance on a non-holonomic space is 
inspired by the work from 1909 of Carath\'eodory \cite{cara} on the mathematical formulation of
thermodynamics.
The modern study of non-holonomic spaces from the viewpoint of 
distance geometry, also known as sub-riemannian geometry, advanced steadily 
due to different lines of research:  hypoelliptic operators
H\"ormander \cite{hormander},  harmonic analysis on homogeneous groups
Folland, Stein \cite{fostein},  probability theory on groups Hazod \cite{hazod}, Siebert \cite{siebert}, studies in geometric
analysis in metric spaces  in relation with
sub-riemannian geometry   Bella\"{\i}che   \cite{bell}, Gromov \cite{gromovsr}, 
 groups with polynomial growth   Gromov \cite{gromovgr}, or Margulis type rigidity
results   Pansu \cite{pansu}. 

Carnot groups are a kind of non-commutative vector spaces, in the sense that 
the group composition (of "vectors") is non-commutative and the multiplication
by scalars is replaced by a contractive automorphism group. The fact that 
 such a  structure appears in relation with non-holonomic spaces (manifolds endowed with
a completely non-integrable distribution) is very non-trivial and deserves an
explanation. 

 Carnot groups seem to appear in this context as a consequence of a 
deformation theory of pseudo-differential operators. The proofs (maybe the most
complete to be found in the big paper of Bella\"{\i}che   \cite{bell}, see also
Margulis, Mostow \cite{marmos1} \cite{marmos2}) involve
very careful development into series and 
 manipulation of iterated brackets of vector fields, along with 
  estimates of the order of remainders, therefore it is not clear if 
 the Carnot group structure appears there by necessity or it is just the
 remnant manifestation of the mathematical tools used to explore the geometry 
 of sub-riemannian spaces. Related to this question is Gromov 
 proposal, made in his paper \cite{gromovsr} , to look at sub-riemannian spaces 
 "from within", that is to deduce the basic results in sub-riemannian geometry 
 in intrinsic manner. 
 Gromov  proposes to use only the Carnot-Carath\'eodory distance
 for this, but it seems that more than the distance is needed 
 in order to get the Carnot group structure of the metric tangent space at a 
 point in a (regular) sub-riemannian manifold. 
 
 A hint that this is true is a result of Siebert \cite{siebert}, which
 essentially characterizes Carnot groups as groups with a contractive group of 
 automorphisms. In the papers \cite{buligadil1} \cite{buligadil2}
 \cite{buligasr} \cite{buligadil3} we achieve this intrinsic description of a 
 sub-riemannian space,  as a particular example of a dilatation structure on a
 metric space. 
 
  In this paper we try to explain that this is a more general phenomenon, 
  namely a   manifestation of interesting emergent operations obtained from iterated
  deformations (or isotopes) of idempotent right quasigroups. The results
  presented here don't belong to distance geometry, but they are inspired by 
  the previous mentioned studies of dilatation structures on metric spaces \cite{buligadil1} 
  \cite{buligadil2} \cite{buligasr}. 
 
 In the language used in this paper a dilatation structure is basically an 
 idempotent right quasigroup endowed with a compatible distance. We show here
 that some of the results obtained previously in metric spaces are in fact 
 true in the realm of uniform idempotent right quasigroups 
 (definition \ref{deftop}). Therefore our question, concerning the apparition of 
 Carnot groups as models of non-commutative tangent spaces in sub-riemannian
 geometry, gets the following interesting answer: uniform idempotent right
 quasigroups offer a minimal concept in order to have a tangent bundle of a 
 space, therefore a decent differential calculus. This minimal concept allow 
 to construct emergent, more complex, algebraic structures, like the one of
 contractible group or of a symmetric space in the sense of Loos \cite{loos}.

\section{Irqs and induced operations}
\label{sec3}

\begin{defi} A right quasigroup is a set $X$ with a binary operation 
$\circ$ such that for each $a, b \in X$ there exists a unique $x \in X$ such that 
$a \, \circ  \, x \, = \, b$. We write the solution of this equation 
$x \, = \, a \, \bullet \, b$. 

A quasigroup is a set $X$ with a binary operation 
$\circ$ such that for each $a, b \in X$ there exist  unique elements  $x, y \in X$ 
such that $a \, \circ  \, x \, = \, b$ and $y \, \circ \, a \, = \, b$. 
We write the solution of the last  equation 
$y \, = \, b \, / \, a$. 

 An idempotent right quasigroup (irq) is a quasigroup with
the property that the operation $\circ$ is idempotent: for any $x \in X$ we have 
$x \, \circ \, x \, = \, x$. 
\label{defquasigroup}
\end{defi}

\begin{rk}
Maybe the most well known example of a right quasigroup comes from 
knot theory. Indeed, in the unpublished correspondence of J.C. Conway and 
G.C. Wraith from 1959, they used the name "wrack" for  a  self-distributive 
right quasigroup generated by a link diagram. Later, Fenn and Rourke
\cite{fennrourke} proposed the name "rack" instead. The correspondence between 
the notation used here and the notation used by Fenn and Rourke is: 
$$x \, \circ \, y \, = y^x$$
 Joyce \cite{joyce} studied and used a particular case of a rack, 
 named "quandle". Quandles are self-distributive idempotent right quasigroups. In the language of Reidemeister
moves, the axioms of a (rack ;  quandle ; irq) correspond respectively to the 
(2,3 ; 1,2,3 ; 1,2) Reidemeister moves. See further remark \ref{remaquandle2} 
for the relation between quandles and contractible groups. 
\label{remaquandle1}
\end{rk}

An equivalent definition for a idempotent right quasigroup (irq) is the
following: $(X, \circ , \bullet )$ is a irq if $X$ is  a  set 
$X$ endowed with two  operations $\circ$ and $\bullet$ which satisfy the following axioms: for any $x , y \in X$  
\begin{enumerate}
\item[(P1)] \hspace{2.cm} $\displaystyle x \, \circ \, \left( x\, \bullet \,  y \right) \, = \, x \, \bullet \, \left( x\, \circ \,  y \right) \, = \, y$
\item[(P2)] \hspace{2.cm} $\displaystyle x \, \circ \, x \, = \, x \, \bullet \, x \,  = \,  x$
\end{enumerate}

We use these operations to define the sum, difference and inverse 
operations of the irq. 

\begin{defi} 
Let $(X, \circ , \bullet )$ be a irq. For 
any $x,u,v \in X$ we define the following operations: 
\begin{enumerate}
\item[(a)] the difference operation is 
$\displaystyle (xuv) \, = \, \left( x \, \circ \, u \right) \, \bullet \, \left(
x \, \circ \, v \right)$. By fixing the first variable $x$ we obtain the difference operation based at $x$: 
$\displaystyle v \, -^{x} \, u \, =  \, dif^{x}(u,v) \, = \, (xuv)$. 
\item[(b)] the sum operation is $\displaystyle )xuv( \, = \, x \, \bullet \left( \left( x \, \circ \, u \right) \, \circ \, 
 v \right)$.  By fixing the first variable $x$ we obtain the sum operation 
 based at $x$: $\displaystyle u \, +^{x} \, v \, =  \, sum^{x}(u,v) \, = \,
 )xuv($.  
\item[(c)] the inverse operation is  $\displaystyle inv(x,u) \, = \, \left( x \,
\circ \, u \right) \, \bullet \,  x $. By fixing the first variable $x$ we obtain the inverse  operator based at $x$: 
$\displaystyle  -^{x} \, u \, =  \, inv^{x} u \, = \, inv(x,u)$. 
\end{enumerate}
 For any $\displaystyle 
k \in \mathbb{Z}^{\circ} = \mathbb{Z} \setminus \left\{ 0 \right\}$  we define the following operations: 
\begin{enumerate}
\item[-] $\displaystyle x \, \circ_{1}\,  u \, = \, x \, \circ \, u$,  $\displaystyle x \, \bullet_{1}\,  u \, = \, x \, \bullet \, u$, 
\item[-]  for any $k > 0$ let 
$\displaystyle  x \, \circ_{k+1}\,  u \, = \, x \, \circ \left(x \circ_{k} \, u \right)$ and 
$\displaystyle  x \, \bullet_{k+1}\,  u \, = \, x \, \bullet \left(x \bullet_{k} \, u \right)$, 
\item[-] for any $k < 0$ let 
$\displaystyle  x \, \circ_{k}\,  u \, = \, x \bullet_{-k} \, u$ and 
$\displaystyle  x \, \bullet_{k}\,  u \, = \, x \circ_{-k} \, u $. 
\end{enumerate}

\label{dplay}
\end{defi}

\begin{rk}
Let $(X, \circ)$, $(Y, \circ)$ be  (right) quasigroups. An isotopy from $X$ to 
$Y$ is a triple of invertible functions $(\alpha, \beta, \gamma)$ defined from 
$X$ to $Y$ with the property that for any $x, y , z \in X$ we have 
$\alpha(x) \, \circ \, \beta(y) \, = \, \gamma(z)$. 

It is then obvious that the operation $\displaystyle +^{x}$ is isotopic with 
$\circ$, with $\alpha(u) \, = \, x \circ u$, $\displaystyle \beta(v) \, = \, v$ and 
$\gamma(w) \, = \, x \circ w$. 
\label{remk1}
\end{rk}

For any $\displaystyle k \in \mathbb{Z}^{\circ}$  the triple 
$\displaystyle (X, \circ_{k} , \bullet_{k} )$ is a irq. We denote  the difference, sum and inverse operations of $\displaystyle (X, \circ_{k} , \bullet_{k} )$  by the same symbols as the ones used for  $(X, \circ , \bullet )$,  with  a subscript "$k$". 

These operations are interesting because they have properties related with 
group operations.

\begin{prop}
In any irq  $(X,\circ, \bullet)$ we have the relations: 
\begin{enumerate}
\item[(a)] $\displaystyle \left( u \, +^{x} \, v \right) \, -^{x} \, u \, = \, v $
\item[(b)] $\displaystyle u \, +^{x} \, \left( v \, -^{x} \, u \right) \, = \, v $
\item[(c)] $\displaystyle v \, -^{x} \, u \, = \, \left(-^{x} u\right)  \, +^{x \circ u} \, v  $
\item[(d)] $\displaystyle -^{x\circ u} \, \left( -^{x} \, u \right) \, = \, u $
\item[(e)] $\displaystyle u \, +^{x} \, \left( v \, +^{x\circ u} \, w \right) \, = \, \left( u \, +^{x} \, v \right) \, +^{x} \, w $
\item[(f)] $\displaystyle  -^{x} \, u \, = \,  x \, -^{x} \, u $
\item[(g)] $\displaystyle  x \, +^{x} \, u \, = \,  u $
\end{enumerate}
\label{pplay}
\end{prop}

\paragraph{Proof.} By remark \ref{remk1} the operation
 $\displaystyle +^{x}$ is isotopic with $\circ$, therefore (a), (b) are just 
  P1 for the operation  $\displaystyle +^{x}$ .  
  
(a) We apply two times P1,  as explained 
further: 
$$\displaystyle \left( u \, +^{x} \, v \right) \, -^{x} \, u \, = \, \left(x \circ u \right) \, \bullet \, \left( 
x \, \circ \, \left( x \, \bullet \, \left( \left( x \circ u \right) \circ \, v
\right)\right)\right) \, = \, $$ 
$$\, = \,  \left(x \circ u \right) \, \bullet \, \left( \left( x \circ u \right) \circ \, v \right) \, = \, v $$

(b) The proof is similar, only the order of application of P1 is reversed, first for "$\displaystyle \left( x \circ u \right) \, \circ \, \left( \left( x \circ u \right) \, \bullet \right.$ ", then for  "$\displaystyle x \, \circ \, \left( x \, \bullet \right.$ " :
$$\displaystyle u \, +^{x} \, \left( v \, -^{x} \, u \right) \, = \, x \, \bullet \, \left( 
\left(x \circ u \right) \, \circ \, \left( \left(x \circ u \right) \, \bullet \, \left( x \, \circ \, v \right)\right)\right) \, = \, 
 x \, \bullet \, \left( x \,  \circ \, v \right) \, = \, v $$

(c) Let us denote by $A$ the expression $\displaystyle A \, = \, \left(x \circ u \right) \, \circ \, \left( -^{x} u \right)$. This expression enters in the right hand side of the equality from (c): 
$\displaystyle  \left(-^{x} u\right)  \, +^{x \circ u} \, v  \, = \, \left(x
\circ u \right) \, \bullet \, \left( A \, \circ \, v \right)$. 
We compute the expression A, using P1:  
$\displaystyle A \, = \, \left(x \circ u \right) \, \circ \, \left( -^{x} u \right) \, = \, \left(x \circ u \right) \, \circ \, \left( \left( x \circ u\right) \bullet \, x \right) \, =  \, x$. 
We use then $A = x$ in the right hand side of the equality (c) and we obtain: 
$$\displaystyle \left(-^{x} u\right)  \, +^{x \circ u} \, v  \, = \, \left(x \circ u \right) \, \bullet \, \left( x \, \circ \, v \right) \, = \, v \, -^{x} \, u $$

(d)  We use again the expression $A$ from the previous computation and the fact that $A = x$, then 
we use P1: 
$\displaystyle  -^{x\circ u} \, \left( -^{x} \, u \right) \, = \, A \, \bullet \, \left( x \circ u\right) \, = \, 
 x \, \bullet \, \left( x \circ u\right) \, = \, u$. 
 
 (e) We compute the left hand side (LHS) and the right hand side (RHS) separately. 
 $$\displaystyle LHS \, = \, x \, \bullet \, \left\{ \left( x \circ u \right) \, \circ \left[   
 \left( x \circ u \right) \bullet \left(  \left( \left( x \circ u \right) \circ \, v \right) \circ \, w \right) \right] \right\} \, = \, 
  x \, \bullet \, \left(  \left( \left( x \circ u \right) \circ \, v \right) \circ \, w \right)$$
 $$RHS \, = \, x \, \bullet \, \left\{  \left[ x \, \circ \, \left( x \, \bullet \, \left( \left(x \circ u \right) \circ \, v      \right)  \right)  \right]  \, \circ \, w  \right\} = \, 
  x \, \bullet \, \left(  \left( \left( x \circ u \right) \circ \, v \right) \circ \, w \right) $$
Therefore $LHS = RHS$. 

(f) Here we use  P2: 
$\displaystyle   x \, -^{x} \, u \, = \, \left( x \circ u \right) \, \bullet  \,
\left( x \circ x \right) \, =$ 
$\displaystyle 
\left( x \circ u \right) \, \bullet \, x \, = \, -^{x} \, u$. 

(g) We use P2, then P1: 
$\displaystyle x \, +^{x} \, u \, = \, x \, \bullet \, \left( \left( x \circ x \right)  \, \circ \, u  \right) \, = \, 
 x \, \bullet \, \left( x  \, \circ \, u  \right) \, = \, u$. 
The proof is done.  \quad $\square$

Let us comment the relations from proposition \ref{pplay}. If we replace all  the superscripts "$x \circ u$" of the signs "$+$" and "$-$" by the superscript "$x$" in the relations (a) to (g)  then we see some interesting patterns. 

For example (e) becomes $\displaystyle u \, +^{x} \, \left( v \, +^{x} \, w \right) \, = \, \left( u \, +^{x} \, v\right) \, +^{x} \, w$, which expresses the associativity of the operation $\displaystyle +^{x} $ . Therefore (e) is a generalized associativity relation. 

Relation (d) takes the form $\displaystyle -^{x} \left( -^{x} u \right) \, = \, u$, that is the inverse 
operator is an involution. 

Relation (c) becomes $\displaystyle v \, -^{x} \, u \, = \, \left(-^{x} u\right)  \, +^{x} \, v  $, which can be seen as a definition of the expression of the left hand side.  Otherwise said $\displaystyle v \, -^{x} \, u $ is the left translate of $v$ (with respect to $\displaystyle +^{x} $ ) by the inverse of $u$. 

Relations (a) and (b) tell us that the inverse of the left translate by $u$ is the left translate by 
$\displaystyle  -^{x} \, u$ (here we don't need to change the superscript "$x$"). 

Relation (g) is transformed into $\displaystyle  x \, +^{x} \, u \, = \,  u $, that is $x$ is a left neutral element 
of $\displaystyle +^{x} $ . 

Without changing the superscripts, let us take $v = x$ in  (c) and use (f): 
$$ \displaystyle  -^{x} \, u \, = \, \left(-^{x} u\right)  \, +^{x \circ u} \, v  $$
If we modify the superscripts then we get that $x$ is also a right neutral element. All in all 
we obtain that $\displaystyle +^{x} $ is  a group operation, with neutral element $x$ and inverse 
 $\displaystyle -^{x} $ . 

Another group of relations will be useful further. 

\begin{prop}
In any irq $(X, \circ , \bullet )$ the following relations are true:
\begin{enumerate}
\item[(h)] $\displaystyle u \, -^{x} \, u \, = \,  (xuu) \, = \, x \, \circ \, u$
\item[(i)]  $\displaystyle u \, -^{x} \, x \, = \,  (xxu) \, =  \, u$
\item[(j)]  $\displaystyle \left( \left(xuu\right) \, \left(xuv\right) \, \left(xuw\right) \right) \, = \, (xvw)$ . 
\item[(k)] Finally, with the notations from definition \ref{dplay},   for any $\displaystyle p, q \in \mathbb{Z}^{\circ}$ and any $x, u , v \in X$ we have the distributivity property: 
$$\displaystyle (x \circ_{q} v) \, -_{p}^{x} \, ( x \circ_{q} u) \, = \, \left( x \circ_{pq} u \right) \, 
\circ_{q} \, \left( v \, -_{pq}^{x} \, u \right)$$

\end{enumerate}
\label{p2play}
\end{prop}

\paragraph{Proof.}
(h) and (i) come from straightforward computations using P2 and  the definition of $(xuv)$. 

(j) Here we use (h) in the left hand side of the equality: 
$$\displaystyle \left( \left(xuu\right) \, \left(xuv\right) \, \left(xuw\right) \right) \, = \,  
\left( \left(x \circ u \right) \, \left(xuv\right) \, \left(xuw\right) \right) \, =
\,  $$ 
$$= \, \left( \left(x \circ u \right) \, \left(\left(x \circ u \right) \bullet \left(x \circ v\right) \right) \, \left(xuw\right) \right) \, = \,  $$
$$ = \, \left( \left( x\circ u\right) \circ       \left(\left(x \circ u \right) \bullet \left(x \circ v\right) \right)                   \right)  \bullet \left(    \left( x\circ u\right) \circ  \left(xuw\right)  \right)   \, = \, $$ 
$$ = \, \left( x \circ v \right) \bullet   \left(    \left( x\circ u\right) \circ  \left(\left( x \circ u \right) \bullet \left( x \circ w \right)\right)  \right) \, = \, \left( x \circ v \right) \bullet \left( x \circ w \right) \, = \, (xvw)$$

 For proving (k) we compute: 
$$\displaystyle \left( x \, ( x \circ_{q} u) \, (x \circ_{q} v) \right)_{p}Ê\, = \, 
\left( x \circ_{p} \, ( x \circ_{q} u) \right)  \, \bullet_{p} \, \left( x \circ_{p} \, (x \circ_{q} v) \right)  \, = \,$$
$$\displaystyle = \, \left( x \circ_{pq} u \right) \, \bullet_{p} \, \left( x \circ_{pq} v \right) \, = \, 
\left( x \circ_{pq} u \right) \, \circ_{q} \, \left(  \left( x \circ_{pq} u \right) \, \bullet_{pq} \, \left( x \circ_{pq} v \right) \right) \, = \, $$ 
 $$ \, = \,  \left( x \circ_{pq} u \right) \, 
\circ_{q} \, \left( x u v \right)_{pq}  \quad  \quad  \square$$

\section{Contractible groups}
\label{sec4}

We can construct irqs from groups, as explained in the following definition. 

\begin{defi}
Let $G$ be a group with operation $(x, y) \mapsto xy$, inverse $\displaystyle x \mapsto x^{-1}$ and neutral element $e$, and $\delta : G \rightarrow G$ a bijective function such that $\delta(e) = e$. We construct the irq $\displaystyle G(\delta) \, = \,  
(G, \circ , \bullet)$ defined by: 
$$ x \, \circ \, u \, = \, x \, \delta \left( x^{-1} u \right) \quad , \quad x \, \bullet \, u \, = \, x \, \delta^{-1} \left( x^{-1} u \right)$$
\end{defi}

For example, if $\displaystyle G = \mathbb{R}^{n}$ with addition and $\varepsilon > 0$, then 
let $\displaystyle \delta_{\varepsilon} (x) \, = \, \varepsilon \, x$. The irq operations are then 
$$ x \, \circ_{k} \, u \, = \, x \, + \, \varepsilon^{k} \left( - x  + u \right) \quad , \quad x \, \bullet_{k} \, u \, = \, x \, + \, \varepsilon^{-k} \left( - x  + u \right)$$
For any $x$ and $u$the point  $x \circ u$ is the result of the homothety of coefficient $\varepsilon$, based at $x$ and applied to $u$. If for example $\displaystyle \varepsilon = \frac{1}{2}$ then $x \circ u$ is the middle point 
between $x$ and $u$, or the arithmetic average of $x$ and $u$. 

The induced operations are: 
$$ v \, -_{k}^{x} \, u \, = \, (x u v )_{k} \, = \, x \, + \, \varepsilon^{k} \left( -x + u \right) \, - \, u \, + \, v $$
$$ u \,+_{k}^{x} \, v \, = \, \,  )x u v (_{k}  \, \, = \, u \, + \, \varepsilon^{k} \left( -u + x \right) \, - \, x \, + \, v $$
If we neglect the term in $\displaystyle \varepsilon^{k}$ then the difference and the sum operations 
are the translates by $x$  of the difference and sum operation in $\displaystyle \mathbb{R}^{n}$. 
This is a general phenomenon which is related with contractible groups, as
explained by proposition \ref{pfirstop}, further on.

In the case of the general irq  $G(\delta)$, if $\delta$ is a group morphism then we get 
 for the induced operations the same expressions as previously: 
$$ v \, -_{k}^{x} \, u \, = \, (x u v )_{k} \, = \, x \,  \delta^{k} \left( x^{-1} u \right) \,   u^{-1} \,  v $$
$$ u \,+_{k}^{x} \, v \, = \, \,  )x u v (_{k}  \, \, = \, u \,  \delta^{k} \left( u^{-1}  x \right) \,  x^{-1}  \, v $$

One particular case is the one of contractible groups.

\begin{defi}
A contractible group is a pair $(G,\alpha)$, where $G$ is a  
topological group with neutral element denoted by $e$, and $\alpha \in Aut(G)$ 
is an automorphism of $G$ such that: 
\begin{enumerate}
\item[-] $\alpha$ is continuous, with continuous inverse, 
\item[-] for any $x \in G$ we have the limit $\displaystyle 
\lim_{n \rightarrow \infty} \alpha^{n}(x) = e$. 
\end{enumerate}
\label{defunu}
\end{defi}

For a contractible group  $(G,\alpha)$, the automorphism $\alpha$ is compactly contractive (Lemma 1.4 (iv) \cite{siebert}), that is: 
for each compact set 
$K \subset G$ and open set $U \subset G$, with $e \in U$, there is 
$\displaystyle N(K,U) \in \mathbb{N}$ such that for any $x \in K$ and $n \in
\mathbb{N}$, $n \geq N(K,U)$, we have $\displaystyle \alpha^{n}(x) \in U$. 

If $G$ is locally compact then $\alpha$ compactly contractive is equivalent with:  each identity neighbourhood 
of $G$ contains an $\alpha$-invariant neighbourhood.

A even more particular case is the one of locally compact groups admitting 
a contractive automorphism group. We begin with  the
definition of a contracting automorphism group \cite{siebert}, definition 5.1. 

\begin{defi}
Let $G$ be a locally compact group. An automorphism group on $G$ is a family 
$\displaystyle T= \left( \tau_{t}\right)_{t >0}$ in $Aut(G)$, such that 
$\displaystyle \tau_{t} \, \tau_{s} = \tau_{ts}$ for all $t,s > 0$. 

The contraction group of $T$ is defined by 
$$C(T) \ = \ \left\{ x \in G \mbox{ : } \lim_{t \rightarrow 0} \tau_{t}(x) = e
\right\} \quad .$$
The automorphism group $T$ is contractive if $C(T) = G$. 
\end{defi}

Next is proposition 5.4 \cite{siebert}, which gives a description of locally
compact groups which admit a contractive automorphism group. 

\begin{prop}
For a locally compact group $G$ the following assertions are equivalent: 
\begin{enumerate}
\item[(i)] $G$ admits a contractive automorphism group;
\item[(ii)] $G$ is a simply connected Lie group whose Lie algebra admits a 
positive graduation.
\end{enumerate}
\label{psiebert}
\end{prop}

The proof of the next proposition is an easy  application of 
the previously explained facts. 

\begin{prop}
Let $(G,\alpha)$ be a contractible group and $G(\alpha)$ be the associated irq. 
Then 
$$\displaystyle \lim_{k \rightarrow \infty} v \, -_{k}^{x} \, u \, = \, 
x \,     u^{-1} \,  v \quad , \quad \lim_{k \rightarrow \infty} u \,+_{k}^{x} \, v \,  = \, u \,   x^{-1}  \, v $$
uniformly with respect to $x, u, v$ in a compact set. 

If $G$ is a group which admits a contractive automorphism group 
$\displaystyle T= \left( \tau_{t}\right)_{t >0}$ then we define for any $t >
0$ the irq  $\displaystyle G(\tau(t))$ which has the associated operations 
denoted respectively by $\displaystyle  v \, -_{t}^{x} \, u$ and $\displaystyle 
u \,+_{k}^{x} \, v$. Then we have: 
$$\displaystyle \lim_{t \rightarrow 0} v \, -_{t}^{x} \, u \, = \, 
x \,     u^{-1} \,  v \quad , \quad 
\lim_{t \rightarrow 0} u \,+_{t}^{x} \, v \,  = \, u \,   x^{-1}  \, v $$
uniformly with respect to $x, u, v$ in a compact set.
\label{pfirstop}
\end{prop}

A particular class of locally compact groups which admit a contractive
automorphism group is  made by Carnot groups. They are related to sub-riemannian or 
Carnot-Carath\'eodory geometry, which is the study of non-holonomic manifolds
endowed with a Carnot-Carath\'eodory distance. Non-holonomic spaces were
discovered in 1926 by  G. Vr\u anceanu \cite{vra1},
\cite{vra2}. The  Carnot-Carath\'eodory distance on a non-holonomic space is 
inspired by Carath\'eodory \cite{cara} work from 1909   on the mathematical formulation of
thermodynamics. Such spaces appear in applications
to thermodynamics, to the mechanics of non-holonomic systems, in the study of
hypo-elliptic operators cf. H\"ormander \cite{hormander}, in harmonic analysis on
homogeneous cones cf. Folland, Stein \cite{fostein}, and
as boundaries of CR-manifolds. We briefly describe the structure of Carnot
groups. 

\begin{defi}
A Carnot (or stratified homogeneous) group is a pair $\displaystyle 
(N, V_{1})$ consisting of a real 
connected simply connected group $N$  with  a distinguished subspace  
$V_{1}$ of  the Lie algebra $Lie(N)$, such that  the following   
direct sum decomposition occurs: 
$$n \ = \ \sum_{i=1}^{m} V_{i} \ , \ \ V_{i+1} \ = \ [V_{1},V_{i}]$$
The number $m$ is the step of the group. The number $\displaystyle Q \ = \ \sum_{i=1}^{m} i 
\ dim V_{i}$ is called the homogeneous dimension of the group. 
\label{dccgroup}
\end{defi}

Because the group is nilpotent and simply connected, the
exponential mapping is a diffeomorphism. We shall identify the 
group with the algebra, if is not locally otherwise stated.

The structure that we obtain is a set $N$ endowed with a Lie
bracket and a group multiplication operation, related by the 
Baker-Campbell-Hausdorff formula. Remark that  the group operation
is polynomial.

Any Carnot group admits a contractive automorphism group.  For any 
$\varepsilon > 0$, the associated morphism is: 
$$ x \ = \ \sum_{i=1}^{m} x_{i} \ \mapsto \ \delta_{\varepsilon} x \
= \ \sum_{i=1}^{m} \varepsilon^{i} x_{i}$$
This is a group morphism and a Lie algebra morphism.

\section{Uniform irqs}
\label{unifirq}

Motivated by proposition \ref{pfirstop}, we introduce  the notion of 
uniform  idempotent right quasigroup.

\begin{defi}
A uniform irq $(X, \circ, \bullet)$ is a separable uniform  
space $X$ endowed with continuous irq operations $\circ$, $\bullet$ such that: 
\begin{enumerate}
\item[-] the operation $\circ$ is compactly contractive: for each compact set 
$K \subset X$ and open set $U \subset X$, with $x \in U$, there is 
$\displaystyle N(K,U) \in \mathbb{N}$ such that for any $u \in K$ and $n \in
\mathbb{N}$, $n \geq N(K,U)$, we have $\displaystyle x \circ_{n} u \in U$; 
\item[-] the following limits exist 
$$ \lim_{k \rightarrow \infty} v \, -_{k}^{x} \, u \, = \, 
v \, -_{\infty}^{x} \, u \quad , \quad \lim_{k \rightarrow \infty} u \,+_{k}^{x}
\, v \,  = \, u \,+_{\infty}^{x} \, v $$
and are uniform with respect to $x, u, v$ in a compact set. 
\end{enumerate}
\label{deftop}
\end{defi}

The main property of a uniform irq is the following. 

\begin{thm}
Let $(X, \circ, \bullet)$ be a uniform irq. Then for any $x \in X$ the operation
$\displaystyle (u,v) \mapsto u \,+_{\infty}^{x} \, v$ gives $X$ the structure of
a contractible group with the contraction $\displaystyle \alpha(u) = x \circ u$.
\label{mainthm}
\end{thm}

\paragraph{Proof.} 
Let us use  the relations from proposition \ref{pplay}, for the 
irq $\displaystyle (X, \circ_{k} , \bullet_{k})$. 
The uniformity assumptions from definition \ref{deftop} 
allow us to pass to the limit with $k \rightarrow \infty$ in these relations. We
obtain the following list: 
\begin{enumerate}
\item[(a)] $\displaystyle \left( u \, +^{x}_{\infty} \, v \right) \,
-^{x}_{\infty} \, u \, = \, v $
\item[(b)] $\displaystyle u \, +^{x}_{\infty} \, \left( v \, -^{x}_{\infty} \, u \right) \, = \, v $
\item[(c)] the inverse $\displaystyle  -^{x}_{\infty} \, u$ (see (f) below) 
exists and
$\displaystyle v \, -^{x}_{\infty} \, u \, = \, \left(-^{x}_{\infty} u\right) 
\, +^{x}_{\infty} \, v  $
\item[(d)] $\displaystyle -^{x}_{\infty} \, \left( -^{x}_{\infty} \, u \right) \, = \, u $
\item[(e)] $\displaystyle u \, +^{x}_{\infty} \, \left( v \, +^{x}_{\infty} \, w
\right) \, = \, \left( u \, +^{x}_{\infty} \, v \right) \, +^{x}_{\infty} \, w $
\item[(f)] $\displaystyle  -^{x}_{\infty} \, u \, = \,  x \, -^{x}_{\infty} \, u
$ exists as a particular case of the limit $\displaystyle  x \, -^{x}_{\infty}
\, u$ with $v = x$, 
\item[(g)] $\displaystyle  x \, +^{x}_{\infty} \, u \, = \,  u $. 
\end{enumerate}
As a consequence $\displaystyle (X , +^{x}_{\infty})$ is a group with 
neutral element $x$ and inverse $\displaystyle u \mapsto \, -^{x}_{\infty} u$.
It is left to prove that $\displaystyle \alpha(u) = x \circ u$ is a group
automorphism. But this is a consequence of passage to the limit in relation 
(k) proposition \ref{p2play}, for $q = 1$ and $p \rightarrow \infty$. \quad
$\square$

\vspace{.5cm}

We have seen that to any uniform irq we can associate a bundle of contractible
groups  $\displaystyle x \in X \mapsto (X, +^{x}_{\infty}, x \circ )$. This
bundle resembles a lot to a tangent bundle, namely: to any $x \in X$ is
associated a contractible group with $x$ as neutral element, which 
 can be seen as the tangent space at $x$. In fact this is a correct picture 
 in the case of a manifold, in the following sense: if we look to a small 
 portion of the manifold then we know that there is a chart of this small 
 portion, which puts it in bijection with an open set in $\mathbb{R}^{n}$. We
 have seen that we can associate to $\mathbb{R}^{n}$ a uniform irq by using 
 as the operation $\circ$ a homothety with fixed ratio $\varepsilon < 1$. This 
 uniform irq is transported on the manifold by the chart. If we ignore the facts
  that we are working not with the whole manifold, but with a small part of it, 
 and not with $\mathbb{R}^{n}$, but with a open set, then indeed we may 
 identify, for any point $x$ in the manifold, a neighbourhood of the point with 
 a neighbourhood of the tangent space at the point, such that the operation 
 of addition of vectors in the tangent space at $x$ is just the operation 
 $\displaystyle +^{x}_{\infty}$ and scalar multiplication by (any integer power 
 of) $\varepsilon$ is just $x \circ$ . 
 
 The same is true in the more complex situation of a sub-riemannian manifold, as
 shown in \cite{buligasr}, in the sense that (locally) we may associate to 
 each point $x$ a "dilatation" of ratio $\varepsilon$, which in turn gives 
 us a structure of uniform irq. In the end we get a bundle of Carnot group 
 operations which can be seen as a tangent bundle of the sub-riemannian
 manifold. (In this case we actually have more structure given by the
 Carnot-Carath\'eodory distance, which induces also a "group norm" on each 
 Carnot group.) 
 
 We shall comment now the fact that a uniform irq can be seen as a
 differential structure. For this it is enough to have a definition of 
 differentiable functions between two uniform irqs. 
 
 \begin{defi}
 Let $(X, \circ, \bullet)$ and $(Y, \circ ,  \bullet \bullet )$ be two 
uniform irqs. A function  $f: X \rightarrow  Y$ is differentiable
 if there is a function $Tf: X \times X \rightarrow Y$ such that 
 $$\lim_{k \rightarrow \infty} f(x) \, \bullet \bullet_{k} \, f (x \circ_{k}
 u) \, = \, Tf (x, u)$$
 uniformly with respect to $x, u$ in compact sets. 
 \label{defder}
 \end{defi}
 This definition corresponds to uniform differentiability in the metric case of 
 dilatation structures (definition 16 and the comments after it in 
 \cite{buligadil1}) and by abstract
 nonsense the application $Tf$ has nice properties, like 
 $\displaystyle Tf(x, \cdot) : (X, +^{x}) \rightarrow (Y, +^{f(x)})$ is a
 morphism of contractible groups.

 It is an interesting question if the emergent operation $\displaystyle 
 +^{x}$ is differentiable. The answer is yes, provided one chooses the 
 right uniform irq structure on $X \times X$, as explained in the metric 
 case in section 8.2 \cite{buligadil1}.  The construction of a uniform irq 
 on $X \times X$, starting from a uniform irq on $X$, such that the 
 operation $\displaystyle  +^{x}$ is differentiable,  is clearly not 
 canonical. But  the operation $\displaystyle +^{x}$ is constructed in a 
 "emergent" manner and this property seems to us stronger than 
 differentiability (although we don't have a proof for this imprecise 
 statement). 
 
 \section{Emergent algebraic structures}
 \label{sec6}

In this section we describe two examples of emergent algebraic structures. 
 We prove in theorem \ref{pgroudlin} 
that contractible groups are in bijection with distributive uniform irqs. 
We then show that such right quasigroups are in fact quasigroups. 
The second example
concerns symmetric spaces in the sense of Loos \cite{loos}. We prove that 
already some properties of the inverse operation in a symmetric space are true
for any uniform irq, then we propose a notion
of uniform symmetric quasigroup, seen as a uniform quasigroup with a
distributivity property for a "approximate inverse" operation $\displaystyle 
\underline{inv}_{k}$, for any $\displaystyle k \in \mathbb{N}^{\circ}$  and prove
that uniform symmetric quasigroups generate a family of symmetric spaces 
which contains the  riemannian symmetric spaces. 

\subsection{Distributive uniform irqs}

 To any group $G$ and contraction  $\alpha$ we can associate a uniform irq, 
as explained by proposition \ref{pfirstop}. The operation $\circ$ is constructed
from $\alpha$ and the group operation. 

We could turn the things upside down and ask if we can recover the contraction
and the group operation from the respective uniform irq. We shall construct the 
initial group operation from the more basic operations of the uniform irq, by a
passage to the limit. In this sense the group operation is "emergent" as the
limit of longer and longer strings of 
compositions of operations $\circ$ and $\bullet$. 

\begin{thm}
Let $(G,\alpha)$ be a contractible group and $G(\alpha)$ be the 
associated uniform irq. Then the irq is distributive, namely it satisfies the
relation 
\begin{equation}
x \circ \left( y \circ z \right) \, = \, \left( x \circ y \right) \circ \left( 
x \circ z \right) 
\label{distributive}
\end{equation}

Conversely, let $(G, \circ, \bullet)$ be a distributive uniform irq. Then there 
is a group operation on $G$ (denoted multiplicatively), with neutral element
$e$, such that: 
\begin{enumerate}
\item[(i)] $\displaystyle xy \, = \, x +^{e}_{\infty} y$ for any $x, y \in G$, 
\item[(ii)] for any $x, y , z \in G$ we have $\displaystyle (xyz)_{\infty} = 
x y^{-1} z$, 
\item[(iii)] for any $x, y \in G$ we have $\displaystyle 
x \circ y \, = \, x (e \circ (x^{-1} y))$. 
\end{enumerate}
In conclusion  there is a bijection between distributive uniform irqs and 
contractible groups. 
\label{pgroudlin}
\end{thm}

\paragraph{Proof.}
Recall that for a contractible group $G, \alpha$ we define 
$\displaystyle x \circ y = x \alpha (x^{-1} y)$. Therefore 
$\alpha(y) = e \circ y$. The distributivity relation (\ref{distributive}) is 
indeed true, as shown by the following string of equalities: 
$$x \circ \left( y \circ z \right) \, = \, x \, \, \alpha (x^{-1} y \, \,  
\alpha (y^{-1} z )) \, = \, x  \, \, \alpha( x^{-1}y) \, \, 
\alpha^{2}(y^{-1} z) \, = \, $$
$$= \, 
 x  \, \, \alpha( x^{-1}y) \,  \, \, \alpha \left( \alpha(y^{-1} x  \, \, x^{-1} z ) \right) \, 
 = x  \, \, \alpha( x^{-1}y) \,  \, \, \alpha \left( \alpha(y^{-1} x) x^{-1}  \,
 \, x \alpha( x^{-1} z ) \right) 
 \, = $$ 
$$= \, (x \circ y) \,  \, \, \alpha \left( \left( x \circ y \right)^{-1} \left(x \circ z
 \right) \right) \, = \, \left( x \circ y \right) \circ \left( 
x \circ z \right) $$  

Conversely, let $e \in X$ and  define 
$\displaystyle xy \, = \, x +^{e}_{\infty} y$. This is a contractible group 
operation, with contraction $\alpha(y) \, e \circ y$. We are left to show 
(ii) and (iii).

(ii). If $(G, \circ , \bullet)$ is a distributive uniform irq then we also have the
following distributivity relations: 
$$x \circ \left( y \bullet z \right) \, = \, \left( x \circ y \right) \bullet \left( 
x \circ z \right)$$ 
$$x \bullet \left( y \circ z \right) \, = \, \left( x \bullet y \right) \circ 
\left( x \bullet z \right)$$
Indeed, the first relation is deduced from (\ref{distributive}) if we replace 
$z$ by $y \bullet z$. The second distributivity relation comes from 
(\ref{distributive}) if we replace $y$ by $x \bullet y$ and $z$ by $x \bullet
z$. In fact we also obtain that for any $k > 0$ 
 $\displaystyle (G, \circ_{k} , \bullet_{k})$ is distributive and moreover, 
 for any $k , l \in \mathbb{Z}$ we have the general distributivity relation 
 $$x \circ_{k} \left( y \circ_{l} z \right) \, = \, \left( x \circ_{k} y \right)
 \circ_{l} \left( x \circ_{k} z \right) $$
 With the help of these relations we obtain the following:  
 $\displaystyle (xyz)_{k} \, = \, x \circ_{k} \left( y \bullet_{k} z \right)$
 and $ \displaystyle \left)xyz\right(_{k} \, = \, y \circ_{k} \left( x
 \bullet_{k} z \right)$. 
 Therefore we get $\displaystyle \left)xyz\right(_{k} \, = \, 
 (yxz)_{k}$. After we pass to the limit with $k \rightarrow \infty$ we obtain: 
 \begin{equation}
 \left)xyz\right(_{\infty} \, = \, (yxz)_{\infty}
 \label{sumdifrel}
 \end{equation}
The following string of equalities is true: 
$$(x e \, (eyz)_{k} \, )_{k} \, = \, x \circ_{k} \left( e \bullet_{k} \left( e
\circ_{k} \left( y \bullet_{k} z \right) \right) \right) \, = \, 
x \circ_{k} \left( y \bullet_{k} z \right) \, = \, (x y z)_{k} $$
After passing to the limit with $k$ we obtain: $ \displaystyle 
(x e \, (eyz)_{\infty} \, )_{\infty} \, = \, (x y z)_{\infty}$, 
therefore we have 
$\displaystyle (x y z)_{\infty} \, = \, (x e \, (eyz)_{\infty} \, )_{\infty} \, = \, 
\left) e x \, (eyz)_{\infty} \, \right(_{\infty} \, = \, x y^{-1} z $. 

(iii).  We start from the following computation: 
$$\left) e \, x \, \left( e \circ \, \left( e x y \right)_{k} \right) \right(_{k}
\, = \, x \circ_{k} \left( e \bullet_{k} \left( e \circ \, \left( e x y
\right)_{k} \right) \right) \, = \, 
x \circ_{k} \left( e \bullet_{k} \left( e \circ \, \left( e \circ_{k} \left(
 x \bullet_{k} y \right) 
\right) \right) \right) \, = $$ 
$$= \, x \circ_{k}  \left( e \circ \left(
 x \bullet_{k} y \right) 
\right) \, = \, \left( x \circ_{k} e \right) \circ \, y $$
We pass to the limit with $k$ and we obtain: 
$$x \circ y \, = \, \left) e \, x \, \left( e \circ \, \left( e x y
\right)_{\infty} \right) \right(_{\infty} \, = \, x (e \circ (x^{-1} y))$$
The proof is done. \quad $\square$

\begin{rk}
In the language of racks and quandles, the previous theorem states that there 
is a bijective correspondence between contracible groups and a family of
quandles. Indeed, as indicated in remark \ref{remaquandle1}, a "distributive 
irq" is another name for a quandle. In this respect, "uniform quandle" could be
an alternative name for a uniform distributive irq. These uniform distributive
irqs are particular examples of topological quandles, studied in \cite{rubin}. 
\label{remaquandle2}
\end{rk}

A loop  is a quasigroup $(X, \circ , \bullet)$ with an identity element $e$
such that $x \circ e \, = \, e \circ x \, = \, x$ for any $x \in X$. Quasigroups are
isotopic with loops, by a well known construction described further: for any 
$x \in X$ we define the operation on $X$ 
\begin{equation}
u \, \circ^{x} \, v \, = \, \left( u \, / \, x \right) \, \circ \, \left( x \, 
\bullet \, v \right)
\label{defisoloop}
\end{equation}
Then $\displaystyle (X, \circ^{x})$ is a loop with identity element $x$. 

Turning back to remark \ref{remk1}, we see that if $X$ is a quasigroup 
the operations $\displaystyle +^{x}$ and $\displaystyle \circ^{x}$ are isotopic. 
 In relation to this, for Carnot groups or more general for normed conical 
 groups (see definition 7 \cite{buligadil2}) we can state the following corollary 
 of proposition 8.4 \cite{buligadil2}. 
 
 \begin{cor}
 Let $G$ be a  contractible group with contraction $\displaystyle
 \delta$ and uniformity coming from a left invariant 
 norm compatible with the contraction.  
   Let us  define 
 the operation $\displaystyle x \, \circ \, y \, = x \, \delta 
 \left( x^{-1} y\ \right)$. Then $\displaystyle (G, \circ_{k})$ is a quasigroup, 
  and moreover, if we denote by  $\displaystyle \circ^{x}_{k}$ the isotope 
  of $\displaystyle \circ_{k}$ as defined by (\ref{defisoloop}), we have  
 $$\lim_{k \rightarrow \infty} u \, \circ_{k}^{x} \, v \, = \, u +^{x} v $$
 uniformly with respect to $x, u , v$ in compact set. 
 \label{corp}
 \end{cor}
 
 \paragraph{Proof.} 
 The fact that  $\displaystyle (G, \circ_{k})$ is a quasigroup is a consequence of 
proposition 8.4 \cite{buligadil2}. Indeed 
the solution of the equation $\displaystyle x \circ_{k} a \, = \, b$ is 
$$x \, = \, b \, /_{k} \, a \, = \,  \prod_{p = 0}^{\infty} \, \delta^{p} 
\left( \delta(
a^{-1}) \, b \right)$$
We can compute then $\displaystyle u \, \circ_{k}^{x} \, v $. We use the
notation $\displaystyle \delta^{x}_{k} y \, = \, x \circ_{k} y$ and we write: 
$$u \, \circ_{k}^{x} \, v \, = \, \left( u \, /_{k} \, x \right) \, \circ_{k} \, 
\left( x \, \bullet_{k} \, v \right) \, = \, \delta_{k}^{u /_{k} x} \, 
\delta_{k^{-1}}^{x} v \, = \, \left( u /_{k} x \right) \, 
\delta_{k} \left( \left( u /_{k} x \right)^{-1} \, x \right) \, x^{-1} v $$
As $k \rightarrow \infty$ we have 
$\displaystyle u/_{k} x \rightarrow u$  uniformly with respect 
to $x, u , v$ in compact set (see for this  the "affine" 
interpretation of $\displaystyle u/_{k} x$ as a ratio  of a collinear triple 
in proposition 8.7  and the estimates provided by 
proposition 8.8 \cite{buligadil2}). We can therefore pass to the limit and 
obtain the conclusion of the corollary. \quad $\square$

\subsection{Uniform symmetric quasigroups}

Further we try to construct symmetric spaces in the sense of Loos \cite{loos} 
from uniform irqs. 

\begin{defi}
$(X, inv)$ is an algebraic symmetric space if $inv : X \times X \rightarrow X$ 
is an operation which satisfies the following axioms: 
\begin{enumerate}
\item[(L1)] $inv$ is idempotent: for any $x \in X$ we have 
$ inv (x,x) = \, x$, 
\item[(L2)] distributivity: for any $x, y , z \in X$ we have 
$$\displaystyle inv (x ,  inv \left( y ,  z \right)) \, = \, 
inv \left( inv \left( x,  y \right) , inv \left( x ,  z \right) \right)$$ 
\item[(L3)] for any $x, y \in X$ we have $inv \left( x , inv 
\left( x  ,  y \right) 
\right) \, = \, y$, 
\item[(L4)] for every $x \in X$ there is a neighbourhood $U(x)$ such that 
 $inv (x ,  y)  \, = \, y$  and $y \in U(x)$ then $x = y$. 
\end{enumerate}
If $X$ is a manifold, $inv$ is smooth (of class $\displaystyle 
\mathcal{C}^{\infty}$) and (L4) is true locally then $(X, inv)$ is a symmetric
space as defined by Loos \cite{loos} chapter II, definition 1. 
\label{defsym}
\end{defi}

We have seen previously that manifolds are particular cases of uniform irqs. 
Our problem is to propose a generalization of Loos symmetric spaces as 
uniform irqs with supplementary algebraic properties. 

For a uniform irq a good candidate for the operation 
$inv$ of an algebraic symmetric space is 
$$ inv_{\infty}(x, y) \, = \, \lim_{k \rightarrow \infty} 
\, inv_{k}(x, y) \quad , \quad inv_{k}(x, y) \, = \, \left( x \circ_{k}  y \right) \bullet_{k} x $$
This limit exists because, as a consequence of proposition \ref{pplay} (f), 
we have   $\displaystyle inv_{\infty}(x, y) \, = \, x -^{x}_{\infty} y$. Indeed,
from corollary \ref{corp} we deduce that axioms (L1), (L4) are true for 
$\displaystyle inv_{\infty}$. Moreover relation (d) proposition \ref{pplay} 
looks like a weak version of axiom (L3). This is confirmed by the following
proposition. 

\begin{prop}
Let $(X, \circ, \bullet)$ be a irq and $T(y, x) \, = \,  (inv(x,y) , x \circ y )$. 
Then $T \circ T \, = \, id$. 
\label{pidem}
\end{prop}

\paragraph{Proof.} 
We define $\delta : X \times X \rightarrow X \times X$ and $i: X \times X
\rightarrow X \times X$ by: 
$$\delta(y,x) \, = \, (x \circ y , x) \quad , \quad i(y,x) \, = \, (x,y)$$
Because $X$ is a irq the transformation $\delta$  is invertible and 
$\displaystyle \delta^{-1}(y,x) \, = \, (x \bullet y , x)$. By direct
computation we have $\displaystyle T \, = \, \delta^{-1} \circ i \circ \delta$. 
The conclusion follows by the obvious fact that $i \circ i \, = \, id$. \quad 
$\square$

There is only one axiom left to investigate, namely (L2). We would like to 
define an operation $\displaystyle \underline{inv}_{k}$, for any 
$\displaystyle k \in \mathbb{N}^{\circ}$, which satisfies (L2) and it 
 is constructed from the operations 
of a uniform irq. 

We have two hints about constructing such an operation. 
The first hint comes from the observation that 
 contractible groups are symmetric spaces and 
 corollary \ref{corp} states that the uniform irq of a contractible group is 
 in fact a quasigroup. The second hint comes from the fact that a Lie group is 
 a symmetric space with operation $\displaystyle inv(u,v)\, = \, u v^{-1}u$. 
 Consider then the uniform irq associated to the Lie group 
 $G$ and a transformation $\displaystyle \delta: G \rightarrow G$, which 
 is not a group morphism, but it is continuous, with continuous inverse and 
 compactly contractive. Then $G(\delta)$ is an uniform irq and after an
 elementary computation we find the interesting relations: 
 $$inv_{k} (u, v) \, = \, inv(u, v \circ_{k} u) \, \quad , \quad inv_{\infty}(u,v)
 \, = \, inv(u,v)$$
These two considerations motivate us to introduce the following definition of 
a uniform symmetric  quasigroup. 

\begin{defi}
A uniform symmetric quasigroup (usq) is a uniform 
irq  $(X, \circ)$ with the following properties: 
\begin{enumerate}
\item[(i)] $\displaystyle (X, \circ_{k}, \bullet_{k})$ is a quasigroup for 
 any $\displaystyle k \in \mathbb{N}^{\circ}$.  For any  $a, b \in X$ 
 we denote by $x \, = \, a \, /_{k} \, b$ the solution of 
 the equation $\displaystyle x \, \circ_{k} \, b \, = \, a$. 
 \item[(ii)] the function $\displaystyle \underline{inv}_{k} : X \times X
\rightarrow X$ defined by 
$$ \displaystyle \underline{inv}_{k}(u,v) \, = \, inv_{k}(u , v \, /_{k} \, u) \, = \, 
\left( u \, \circ_{k} \, \left( v \,  /_{k} \,  u \right) \right) \, \bullet_{k}
\, u$$
 satisfies (L2) (is distributive). 
\item[(iii)] the function $\displaystyle inv_{\infty} $ satisfies (L4).
\end{enumerate}
\label{defsymirq}
\end{defi}

In order to cover all symmetric spaces 
(like the compact ones for example) we should ask that condition (i) is true 
only locally. (The locality condition in (iii) is of a different nature, see the
following remark). In this paper we  make only  local considerations.  The problem of locality in (i)  can be solved, but this 
leads to more involved reasoning and notations, which only obscure the 
ideas presented here. For a way to treat the locality problem (but within the
metric space theory) see \cite{buligadil1} \cite{buligadil2}. 

\begin{rk}
The condition (iii) is not so severe. For example if the tangent space at 
$x$, namely $\displaystyle (X, +^{x})$, 
is a Carnot group (which is the case for the symmetric spaces in the sense of 
Loos) then (iii) is trivially true. 

Thus a more restrictive condition that   (iii) is: for every $x \in X$ the 
group $\displaystyle (X, +^{x})$ is locally compact (with respect to the
topology inherited from the uniformity on $X$)  and admits a contractive 
automorphism group $C(x)$ such  that the application 
$\displaystyle \delta^{x}(y) \, = \, x \, \circ \, y$ belongs to $C(x)$. 
\label{rkl4}
\end{rk}

\begin{prop}
If $(X, \circ, \bullet)$ is  a uniform symmetric quasigroup then $\displaystyle 
(X, inv_{\infty})$ is an algebraic symmetric space. Conversely, let $X$ be 
any riemannian symmetric space, $\varepsilon \in (0,1)$, 
let the operation $\circ$ be defined by   
$\displaystyle x \, \circ \, exp^{x}(X) \, = \, exp^{x}(\varepsilon \, X)$, where 
$exp$ is the geodesic exponential. Then $X$ is a uniform symmetric quasigroup. 
\label{thmsym}
\end{prop}

\paragraph{Proof.} With $\displaystyle inv \, = \, inv_{\infty}$ the properties 
(L1), (L3) are true for any usq. (L4) is just the condition 
(iii). It is left to study the distributivity (L2). This is straightforward
because $\displaystyle \underline{inv}_{k}(x, u \circ_{k} x) \, = \, inv_{k}(x, u)$,
and $\displaystyle inv_{k}(x, \cdot)$ converges uniformly on compact sets to 
$\displaystyle inv(x, \cdot)$,  therefore we have the uniform limit 
$$\lim_{k \rightarrow \infty} \underline{inv}_{k}(x, v) \, = \, inv(x,v)$$
But then we can pass to the limit in the distributivity relation 
for $\displaystyle \underline{inv}_{k}$ and eventually obtain the 
distributivity condition (L4) for $\displaystyle inv_{\infty}$. 

For the second part of the proposition we first recall that we are 
making only local considerations, therefore we shall consider that 
any two points are joined by a unique geodesic. The fact that 
$(X, \circ)$ is a uniform irq reduces to an elementary computation using the
geodesic exponential) (basically that the derivative of 
geodesic exponential $\displaystyle 
exp^{x}(\varepsilon X)$ with respect to $\varepsilon$ at $\varepsilon = 0 $
equals $X$).  We only have to check then  if the operation
$\circ$ is a quasigroup operation and if it satisfies (ii) \ref{defsymirq}. But this
is obvious because we are using the geodesic exponential for the definition of
this operation. Then $x$, $x \circ u$ and $u$ are on the same geodesic and with
respect to the coordinate system with origin in $x$ along this geodesic, 
the coordinate of $x \circ u$ equals $\varepsilon$ times the coordinate of $y$. This
shows that $x$ is determined by $x \circ u$ and $u$. Concerning the condition 
(ii), we claim that $\displaystyle \underline{inv}_{k}(u,v) \, = \, inv(u,v)$
for any $\displaystyle k \in \mathbb{N}^{\circ}$. This is equivalent with 
$\displaystyle inv_{k}(u, v) \, = \, inv(v \, \circ_{k} \, u)$. 
Indeed let $\gamma$ be the geodesic joining $u$ with $v$, parameterized by 
length, with origin in $u$. If $v$ has coordinate $V$, then $ \displaystyle 
u \, \circ_{k} \, v$ is on $\gamma$, with coordinate  
 coordinate $\displaystyle  \varepsilon^{k} V$ and $\displaystyle inv_{k}(u, v)$
 is then also on $\gamma$, with coordinate $\displaystyle (\varepsilon^{k} - 1)
 V$. But $\displaystyle inv(v \, \circ_{k} \, u)$ is on $\gamma$ as well, with the 
 same coordinate as $\displaystyle inv_{k}(u, v)$. \quad $\square$

\end{document}